\documentclass{amsart}
\usepackage{amsfonts,amsmath,amsthm,amssymb}
\newtheorem{theorem}{Theorem}
\newtheorem{corollary}{Corollary}
\newtheorem{proposition}{Proposition}
\title[Integral group ring of the Mathieu simple group $M_{12}$]
{Integral group ring of the \\ Mathieu simple group $M_{12}$}
\date{}

\author{V.A.~Bovdi, A.B.~Konovalov and S.~Siciliano}

\address{V.A.~Bovdi
\newline Institute of Mathematics, University of Debrecen,
\newline P.O. Box 12, H-4010 Debrecen, Hungary
\newline Institute of Mathematics and Informatics, College of Ny\'\i regyh\'aza,
\newline S\'ost\'oi \'ut 31/b, H-4410 Ny\'\i regyh\'aza, Hungary}
\email{vbovdi@math.klte.hu}

\address{A.B.~Konovalov
\newline School of Computer Science, University of St Andrews
\newline Jack Cole Building, North Haugh, St Andrews, Fife, KY16 9SX, Scotland}
\email{konovalov@member.ams.org}

\address{S.~Siciliano
\newline Dipartimento di Matematica "E. De Giorgi", 
         Universit\`{a} degli Studi di Lecce,
\newline Via Provinciale Lecce-Arnesano, 73100-LECCE, Italy}
\email{salvatore.siciliano@unile.it}

\thanks{The research was supported by OTKA grants No.T 43034,
No.K61007 and Francqui Stichting (Belgium) grant ADSI107}

\subjclass{Primary 16S34, 20C05, secondary 20D08}

\thanks{}

\keywords{Zassenhaus conjecture, Kimmerle conjecture,
torsion unit, partial augmentation, integral group ring}

\begin{document}
\begin{abstract}
We consider the  Zassenhaus conjecture  for the normalized unit
group of the integral group ring of  the Mathieu sporadic group
$M_{12}$. As a consequence, we confirm for this group the 
Kimmerle's conjecture on prime graphs.
\end{abstract}

\maketitle

\section{Introduction, conjectures   and main results}
\label{Intro}

Let $V(\mathbb Z G)$ be  the normalized unit group of the
integral group ring $\mathbb Z G$ of  a finite group $G$. A long-standing
conjecture  of H.~Zassenhaus {\bf (ZC)} says that every torsion unit
$u\in V(\mathbb ZG)$ is conjugate within the rational group algebra
$\mathbb Q G$ to an element in $G$.

For finite simple groups the main tool for the investigation of the
Zassenhaus conjecture is the Luthar-Passi method, introduced in
\cite{Luthar-Passi} to solve it for $A_{5}$. 
Later M.~Hertweck in \cite{Hertweck1} extended
the Luthar-Passi method and applied it for the investigation
of the Zassenhaus conjecture for $PSL(2,p^{n})$.
The Luthar-Passi method proved to be useful for groups containing 
non-trivial normal subgroups as well. For some recent results we refer to
\cite{Bovdi-Hofert-Kimmerle,Bovdi-Konovalov,
Hertweck2,Hertweck1,Hertweck3,Hofert-Kimmerle}. 
Also, some related properties and  some weakened variations of the
Zassenhaus  conjecture can be found in
\cite{Artamonov-Bovdi,Luthar-Trama} and
\cite{Bleher-Kimmerle,Kimmerle}.

First of all, we need to introduce some notation. By $\# (G)$ we
denote the set of all primes dividing the order of $G$. The
Gruenberg-Kegel graph (or the prime graph) of $G$ is the graph
$\pi (G)$ with vertices labeled by the primes in $\# (G)$ and with
an edge from $p$ to $q$ if there is an element of order $pq$ in
the group $G$. In \cite{Kimmerle} W.~Kimmerle   proposed the
following weakened variation of the Zassenhaus conjecture:

\begin{itemize}
\item[]{\bf (KC)} \qquad
If $G$ is a finite group then $\pi (G) =\pi (V(\mathbb Z G))$.
\end{itemize}

In particular, in the same  paper  W.~Kimmerle verified   that
{\bf (KC)} holds for finite Frobenius and solvable groups. Note
that with respect to the so-called $p$-version of the Zassenhaus 
conjecture the investigation
of Frobenius groups was completed by  M.~Hertweck and the first
author  in \cite{Bovdi-Hertweck}. In
\cite{Bovdi-Jespers-Konovalov, Bovdi-Konovalov, Bovdi-Konovalov-M23} 
{\bf (KC)} was confirmed for sporadic simple groups $M_{11}$, $M_{23}$ and 
some Janko simple groups.

Here  we continue these investigations for the Mathieu simple group
$M_{12}$. Although using the Luthar-Passi method we cannot prove the
rational conjugacy for torsion units of $V(\mathbb Z M_{12})$, our
main result gives a lot of information on partial augmentations of
these units. In particular, we confirm the Kimmerle's conjecture for
this group.

Let $G=M_{12}$. It is well known
(see \cite{GAP,Gorenstein}) that $|G|=2^6\cdot 3^3\cdot 5\cdot 11$ and
$exp(G)=2^3 \cdot 3 \cdot  5 \cdot 11$.
Let
$$
\mathcal{C} =\{ C_{1}, C_{2a} , C_{2b} , C_{3a} , C_{3b} , C_{4a}
, C_{4b} , C_{5a}, C_{6a} , C_{6b} , C_{8a} , C_{8b} , C_{10a} ,
C_{11a} , C_{11b}\}
$$
be the collection of all conjugacy classes of $M_{12}$, where the first
index denotes the order of the elements of this conjugacy class
and $C_{1}=\{ 1\}$. Suppose $u=\sum \alpha_g g \in V(\mathbb Z G)$
has finite order $k$. Denote by
$\nu_{nt}=\nu_{nt}(u)=\varepsilon_{C_{nt}}(u)=\sum_{g\in C_{nt}}
\alpha_{g}$ the partial augmentation of $u$ with respect to
$C_{nt}$. From S.D.~Berman's Theorem \cite{Berman} one knows that
$\nu_1 =\alpha_{1}=0$ and
\begin{equation}\label{E:1}
\sum_{C_{nt}\in \mathcal{C}} \nu_{nt}=1.
\end{equation}
Hence, for any character $\chi$ of $G$, we get that $\chi(u)=\sum
\nu_{nt}\chi(h_{nt})$, where $h_{nt}$ is a representative of the
conjugacy class $ C_{nt}$.

Our main result is the following
\begin{theorem}\label{T:1}
Let $G$ denote the Mathieu simple group $M_{12}$.  Let  $u$ be a
torsion unit of $V(\mathbb ZG)$ of order $|u|$. The following
properties hold:

\begin{itemize}
\item[(i)] If  $|u|\not\in\{12,20,24,40\}$, then $|u|$ coincides  with the order
of some element  $g\in G$.
\item[(ii)] If $|u|=2$, then the tuple of the  partial
augmentations of $u$ belongs to the set
\[
\begin{split}
\big\{\;  (\nu_{2a} , \nu_{2b} , \nu_{3a} , \nu_{3b} , \nu_{4a} ,
\nu_{4b} , \nu_{5a}, \nu_{6a} , \nu_{6b} , \nu_{8a} , \nu_{8b} ,
\nu_{10a} , \nu_{11a} , \nu_{11b}) \in \mathbb Z^{14} \mid& \\
(\nu_{2a},\nu_{2b}) \in \{ \;  ( 0, 1 ),\; ( -2, 3 ),\; ( 2, -1 ),\; ( 1,
0 ),\; ( 3, -2 ),\; ( -1, 2 ) \; \},& \\
 \nu_{kx}=0,\;   kx\not\in\{2a,2b\}&\; \big\} .
\end{split}
\]
\item[(iii)]
 If $|u|=3$, then the tuple of the partial augmentations
of $u$ belongs to the set
\[
\begin{split}
\big\{\; (\nu_{2a} , \nu_{2b} , \nu_{3a} , \nu_{3b} , \nu_{4a} ,
\nu_{4b} , \nu_{5a}, \nu_{6a} , \nu_{6b} , \nu_{8a} , \nu_{8b} ,
\nu_{10a} , \nu_{11a} , &\nu_{11b}) \in \mathbb Z^{14} \mid \\
(\nu_{3a},\nu_{3b}) \in \{ \; ( 0, 1 ),\; ( 2, -1 ),\; ( 1, 0 ),\; ( 3,
-2 ),\; ( -1, 2 ) & \; \},\quad  \nu_{kx}=0,\\
& kx\not\in\{3a,3b\}\; \big\} .
\end{split}
\]
\item[(iv)] If $|u|=5$, then $u$ is rationally
conjugate to some  $g\in G$;
\item[(v)] If $|u|=10$, then the tuple of the  partial
augmentations of $u$ belongs to the set
\[
\begin{split}
\big\{\;  (\nu_{2a} , \nu_{2b} , \nu_{3a} , \nu_{3b} , \nu_{4a} ,
\nu_{4b} , \nu_{5a}, \nu_{6a} , \nu_{6b} , \nu_{8a} , \nu_{8b} ,
\nu_{10a} , \nu_{11a} , \nu_{11b}) \in \mathbb Z^{14} \mid& \\
 (\nu_{2a},\nu_{2b},\nu_{10a}) \in \{ \;  (0, 0, 1 ),\; ( 1, 1,
-1)\},\;\quad  \nu_{kx}=0,& \\
kx\not\in\{2a,2b,10a\}& \; \big\}.
\end{split}
\]
\item[(vi)] If $|u|=11$, the tuple of the  partial
augmentations of $u$ belongs to the set
\[
\begin{split}
\big\{\;  (\nu_{2a} , \nu_{2b} , \nu_{3a} , \nu_{3b} , \nu_{4a} ,
\nu_{4b} , \nu_{5a}, \nu_{6a} , \nu_{6b} , \nu_{8a} , \nu_{8b} ,
\nu_{10a} , \nu_{11a} , \nu_{11b}) \in \mathbb Z^{14} \mid &\\
(\nu_{11a},\nu_{11b}) \in \{ \;  ( 0, 1 ),\; ( 2, -1 ),\; ( 1,
0),\; (-1, 2)\;\},\quad \nu_{kx}=0,& \\
   kx\not\in\{11a,11b\}& \; \big\}.
\end{split}
\]
\end{itemize}
\end{theorem}

As an immediate consequence of  part (i) of the Theorem we obtain

\begin{corollary} If $G=M_{12}$ then
$\pi(G)=\pi(V(\mathbb ZG))$.
\end{corollary}

\section{Preliminaries}
The following result is a reformulation of
the Zassenhaus conjecture in terms of partial augmentations
of torsion units.

\begin{proposition}\label{P:5}
(see \cite{Luthar-Passi} and
Theorem 2.5 in \cite{Marciniak-Ritter-Sehgal-Weiss})
Let $u\in V(\mathbb Z G)$
be of order $k$. Then $u$ is conjugate in $\mathbb
QG$ to an element $g \in G$ if and only if for
each $d$ dividing $k$ there is precisely one
conjugacy class $C$ with partial augmentation
$\varepsilon_{C}(u^d) \neq 0 $.
\end{proposition}

The next two results now yield that several partial augmentations
are zero.

\begin{proposition}\label{P:3}(see \cite{Luthar-Passi} and
Theorem 2.7 in \cite{Marciniak-Ritter-Sehgal-Weiss})
Let $u$ be a torsion unit of $V(\mathbb ZG)$. Let
$C$ be a conjugacy class of $G$. If $a\in C$ and
$p$ is a prime dividing  the order of $a$ but not
the order of $u$ then  $\varepsilon_C(u)=0$.
\end{proposition}

\begin{proposition}\label{P:4}
(see \cite{Hertweck2}, Proposition 3.1;
\cite{Hertweck1}, Proposition 2.2)
Let $G$ be a finite
group and let $u$ be a torsion unit in $V(\mathbb
ZG)$. If $x$ is an element of $G$ whose $p$-part,
for some prime $p$, has order strictly greater
than the order of the $p$-part of $u$, then
$\varepsilon_x(u)=0$.
\end{proposition}

Another important restriction on partial augmentations is given by
the next result, explained in details in \cite{Luthar-Passi} and
\cite{Bovdi-Hertweck,Hertweck1}.
\begin{proposition}\label{P:1}
(see \cite{Luthar-Passi,Hertweck1}) Let either $p=0$ or $p$ a prime
divisor of $|G|$. Suppose
that $u\in V( \mathbb Z G) $ has finite order $k$ and assume $k$ and
$p$ are coprime in case $p\neq 0$. If $z$ is a primitive $k$-th root
of unity and $\chi$ is either a classical character or a $p$-Brauer
character of $G$, then for every integer $l$ the number
\begin{equation}\label{E:2}
\mu_l(u,\chi, p ) =
\textstyle\frac{1}{k} \sum_{d|k}Tr_{ \mathbb Q (z^d)/ \mathbb Q }
\{\chi(u^d)z^{-dl}\}
\end{equation}
is a non-negative integer.
\end{proposition}

Note that if $p=0$, we will use the notation $\mu_l(u,\chi
, * )$ for $\mu_l(u,\chi , 0)$.

Finally, we shall use the well-known bound for
orders of torsion units.

\begin{proposition}\label{P:2}  (see  \cite{Cohn-Livingstone})
The order of a torsion element $u\in V(\mathbb ZG)$
is a divisor of the exponent of $G$.
\end{proposition}


\section{Proof of the Theorem}

Throughout this section we denote $M_{12}$ by $G$. 
The character table of $G$,
as well as the $p$-Brauer character tables, which will be denoted by
$\mathfrak{BCT}{(p)}$ where $p\in\{2,3,5,11\}$, can be found using
the computational algebra system GAP \cite{GAP}. For the characters
and conjugacy classes we will use throughout the paper the same
notation, indexation inclusive, as used in GAP.

Since the group $G$ possesses elements of
orders $2$, $3$, $4$, $5$, $6$, $8$, $10$ and  $11$, first of all we
investigate  units of some of these orders (except the units of
orders $4,6$ and $8$). After this, by Proposition \ref{P:2}, the
order of each torsion unit divides the exponent of $G$, so it
remains to consider units of orders $12$, $15$, $20$, $22$, $33$ and
$55$. We prove that no units of all these orders, except for $12$
and $20$, do  appear in $V(\mathbb ZG)$.

Now we consider each case separately.

\noindent$\bullet$ Let $u$ be an involution. By (\ref{E:1}) and
Proposition \ref{P:3} we have that $\nu_{2a}+\nu_{2b}=1$. Applying
Proposition \ref{P:1} to the character $\chi_{2}$,
we get the following system
\[
\begin{split}
\mu_{0}(u,\chi_{2},*) & = \textstyle \frac{1}{2} (-\;t_1 +11) \geq 0; \qquad
\mu_{1}(u,\chi_{2},*)  = \textstyle \frac{1}{2} ( \;t_1 +11) \geq 0; \\ 
\mu_{0}(u,\chi_{2},3) & = \textstyle \frac{1}{2} (-2t_2 +10) \geq 0; \qquad
\mu_{1}(u,\chi_{2},3)  = \textstyle \frac{1}{2} (2t_2 +10) \geq 0, \\ 
\end{split}
\]
where $t_1=\nu_{2a} -3 \nu_{2b}$ and $t_2=\nu_{2a} - \nu_{2b}$.
Obviously, $t_1\in \{2s+1\mid -6\leq s\leq 5\}$ and $-5 \leq t_2\leq
5$. From these restrictions and the requirement that all
$\mu_i(u,\chi_{j},p)$ must be non-negative integers we obtain six
pairs $(\nu_{2a},\nu_{2b})$ listed in part (ii) of Theorem~\ref{T:1}. 
Note that checking conditions of Proposition \ref{P:1}
for all other combinations of $\chi_j$ and $p \in {0,2,3,5,11}$, 
we will not get further restrictions on partial augmentations.

\noindent $\bullet$ Let $u$ be a unit of order $3$. By (\ref{E:1})
and Proposition \ref{P:3} we get $\nu_{3a}+\nu_{3b}=1$. By
(\ref{E:2}) we obtain  the system of inequalities
$$
\mu_{0}(u,\chi_{2},*) = \textstyle \frac{1}{3} (2t_1 +11) \geq 0; \quad
\mu_{0}(u,\chi_{4},*) = \textstyle \frac{1}{3} (-2t_1 +16) \geq 0,
$$
where $t_1=2\nu_{3a} -\nu_{3b}$. Obviously, $t_1\in\{-4,-1,2,
5,8\}$. Using $t_2=\nu_{3a} +\nu_{3b}=1$ and the condition that
$\mu_i(u,\chi_{j},p)$ are non-negative integers, we obtain the five
pairs $(\nu_{3a},\nu_{3b})$ listed in part (iii) of Theorem
\ref{T:1}.

\noindent $\bullet$ Let $u$ be a unit of order $5$. Using
Propositions \ref{P:3} and \ref{P:4} we obtain that all partial
augmentations except one are zero. Thus by Proposition \ref{P:3} the
proof of part (iv) of  Theorem \ref{T:1} is done.

\noindent$\bullet$ Let $u$ be a unit of order $10$. By (\ref{E:1})
and Proposition \ref{P:3} we have that
$$
\nu_{2a}+\nu_{2b}+\nu_{5a}+\nu_{10a}=1.
$$
Now we need to consider the six cases defined by part (ii) of
Theorem \ref{T:1}.

\noindent
Case 1. $\chi(u^5)=\chi(2a)$.
Applying Proposition \ref{P:1}, we get the system of inequalities
\[
\begin{split}
\mu_{0}(u,\chi_{2},*) & = \textstyle \frac{1}{10} (-4 t_1 + 14) \geq 0; \qquad 
\mu_{5}(u,\chi_{2},*)   = \textstyle \frac{1}{10} (4 t_1 + 16) \geq 0; \\ 
\mu_{0}(u,\chi_{4},*) & = \textstyle \frac{1}{10} (4 t_2 + 24) \geq 0; \qquad 
\mu_{5}(u,\chi_{4},*)   = \textstyle \frac{1}{10} (-4 t_2 + 16) \geq 0; \\ 
\mu_{0}(u,\chi_{7},*) & = \textstyle \frac{1}{10} (4 t_3 + 56) \geq 0; \qquad 
\mu_{5}(u,\chi_{7},*)   = \textstyle \frac{1}{10} (-4 t_3 + 44) \geq 0, \\ 
\end{split}
\]
where
$t_1=\nu_{2a}-3\nu_{2b}-\nu_{5a}+\nu_{10a}$,
$t_2=4 \nu_{2a} + \nu_{5a} - \nu_{10a}$
and
$t_3=6\nu_{2a}+6\nu_{2b}-\nu_{5a}+\nu_{10a}$.
Clearly,
$t_1 \in \{-4,1\}$,
$t_2 \in \{-6,-1,4\}$
and
$t_3 \in \{ -14, -9, -4, 1, 6, 11 \}$.
From these restrictions and the condition for
$\mu_i(u,\chi_{j},p)$ to be non-negative integers, we obtain only
three solutions $(-1, 0, 0, 2)$, $(0, 0, 0, 1)$ and $(1, 0, 0, 0)$,
and because of the additional inequality
$\mu_{1}(u,\chi_{4},*) = \textstyle \frac{1}{10} (4 \nu_{2a} +  \nu_{5a} -  \nu_{10a} + 11) \geq 0$ 
there  remains only one solution $(0, 0, 0, 1)$.

\noindent
Case 2. $\chi(u^5)=\chi(2b)$.
Again by using Proposition \ref{P:1}, we obtain
\[
\begin{split}
\mu_{0}(u,\chi_{2},*) & = \textstyle \frac{1}{10} (-4 t_1 + 18) \geq 0; \quad 
\mu_{5}(u,\chi_{2},*)   = \textstyle \frac{1}{10} (4 t_1 + 12) \geq 0; \\ 
\mu_{0}(u,\chi_{4},*) & = \textstyle \frac{1}{10} (4 t_2 + 20) \geq 0; \qquad 
\mu_{5}(u,\chi_{4},*)   = \textstyle \frac{1}{10} (-4 t_2 + 20) \geq 0; \\ 
\mu_{0}(u,\chi_{7},*) & = \textstyle \frac{1}{10} (4 t_3 + 56) \geq 0; \qquad 
\mu_{5}(u,\chi_{7},*)   = \textstyle \frac{1}{10} (-4 t_3 + 44) \geq 0, \\ 
\end{split}
\]
where $t_1$, $t_2$ and $t_3$ are defined as in the previous case.
From this, it follows that
$t_1 \in \{-3,2\}$,
$t_2 \in \{-5,0,5\}$
and
$t_3 \in \{ -14, -9, -4, 1, 6, 11 \}$.
Only three 4-tuples
$(-1, 1, 0, 1)$, $(0, 1, 0, 0)$ and $(1, 1, 0, -1)$ may satisfy
these restrictions and the condition for
$\mu_i(u,\chi_{j},p)$ to be non-negative integers. After
considering additional
inequalities
\[
\begin{split}
\mu_{1}(u,\chi_{4},*) & = \textstyle \frac{1}{10} (4 \nu_{2a} +  \nu_{5a} -  \nu_{10a} + 15) \geq 0; \\ 
\mu_{0}(u,\chi_{4},3) & = \textstyle \frac{1}{10} (12 \nu_{2a} - 4 \nu_{2b} - 8 \nu_{10a} + 14) \geq 0, \\ 
\end{split}
\]
we can eliminate two more solutions, and so there  remains only $(1,
1, 0, -1)$.

\noindent Case 3. $\chi(u^5)=-2\chi(2a)+3\chi(2b)$. As above, by
Proposition \ref{P:1} we obtain that
\[
\begin{split}
\mu_{1}(u,\chi_{2},*) & = \textstyle \frac{1}{10} (- t_1 - 1) \geq 0; \qquad 
\mu_{5}(u,\chi_{2},*)   = \textstyle \frac{1}{10} (t_1 + 4) \geq 0; \\ 
\mu_{0}(u,\chi_{4},*) & = \textstyle \frac{1}{10} (4 t_2 + 12) \geq 0; \qquad 
\mu_{2}(u,\chi_{4},*)   = \textstyle \frac{1}{10} (- t_2 + 7) \geq 0; \\ 
\mu_{0}(u,\chi_{7},*) & = \textstyle \frac{1}{10} (4 t_3 + 56) \geq 0; \qquad 
\mu_{5}(u,\chi_{7},*)   = \textstyle \frac{1}{10} (-4 t_3 + 44) \geq 0, \\ 
\end{split}
\]
where $t_1$, $t_2$ and $t_3$ are defined as in the previous case.
This yields
$t_1=-1$,
$t_2 \in \{-3,7\}$
and
$t_3 \in \{ -14, -9, -4, 1, 6, 11 \}$,
and there are no solutions satisfying these
restrictions and the condition for
$\mu_i(u,\chi_{j},p)$ to be non-negative integers.

\noindent
Case 4. $\chi(u^5)=2\chi(2a)-\chi(2b)$. Again, for the same
$t_1$, $t_2$ and $t_3$ we have
\[
\begin{split}
\mu_{0}(u,\chi_{2},*) & = \textstyle \frac{1}{10} (-4 t_1 + 10) \geq 0; \quad 
\mu_{2}(u,\chi_{2},*)   = \textstyle \frac{1}{10} (t_1 + 5) \geq 0; \\ 
\mu_{0}(u,\chi_{4},*) & = \textstyle \frac{1}{10} (4 t_2 + 28) \geq 0; \qquad 
\mu_{5}(u,\chi_{4},*)   = \textstyle \frac{1}{10} (-4 t_2 + 12) \geq 0; \\ 
\mu_{0}(u,\chi_{7},*) & = \textstyle \frac{1}{10} (4 t_3 + 56) \geq 0; \qquad 
\mu_{5}(u,\chi_{7},*)   = \textstyle \frac{1}{10} (-4 t_3 + 44) \geq 0, \\ 
\end{split}
\]
We have that $t_1=-5$, $t_2 \in \{-7,-2,-3\}$ and $t_3 \in \{ -14,
-9, -4, 1, 6, 11 \}$, and, by the same arguments as in the previous
case, we have no solutions.

\noindent
Case 5. $\chi(u^5)=3\chi(2a)-2\chi(2b)$. Again, we obtain that
\[
\begin{split}
\mu_{0}(u,\chi_{2},*) & = \textstyle \frac{1}{10} (-4 t_1 + 6) \geq 0; \qquad 
\mu_{2}(u,\chi_{2},*)   = \textstyle \frac{1}{10} ( t_1 + 1) \geq 0; \\ 
\mu_{1}(u,\chi_{4},*) & = \textstyle \frac{1}{10} ( t_2 + 3) \geq 0; \qquad\quad 
\mu_{5}(u,\chi_{4},*)   = \textstyle \frac{1}{10} (-4 t_2 + 8) \geq 0; \\ 
\mu_{0}(u,\chi_{7},*) & = \textstyle \frac{1}{10} (4 t_3 + 56) \geq 0; \qquad 
\mu_{5}(u,\chi_{7},*)   = \textstyle \frac{1}{10} (-4 t_3 + 44) \geq 0, \\ 
\end{split}
\]
so
$t_1=-1$,
$t_2=-3$
and
$t_3 \in \{ -14, -9, -4, 1, 6, 11 \}$, which yields no solutions by the same arguments.

\noindent
Case 6. $\chi(u^5)=-\chi(2a)+2\chi(2b)$.
Similarly, we get
\[
\begin{split}
\mu_{1}(u,\chi_{2},*) & = \textstyle \frac{1}{10} (- t_1 + 3) \geq 0; \qquad 
\mu_{5}(u,\chi_{2},*)   = \textstyle \frac{1}{10} (4 t_1 + 8) \geq 0; \\ 
\mu_{0}(u,\chi_{4},*) & = \textstyle \frac{1}{10} (4 t_2 + 16) \geq 0; \qquad 
\mu_{5}(u,\chi_{4},*)   = \textstyle \frac{1}{10} (-4 t_2 + 24) \geq 0; \\ 
\mu_{0}(u,\chi_{7},*) & = \textstyle \frac{1}{10} (4 t_3 + 56) \geq 0; \qquad 
\mu_{5}(u,\chi_{7},*)   = \textstyle \frac{1}{10} (-4 t_3 + 44) \geq 0; \\ 
\end{split}
\]
so $t_1=3$, $t_2 \in \{-4,1,6\}$ and $t_3 \in \{ -14, -9, -4, 1, 6,
11 \}$. As before, in this case too we have no solutions. Thus, part
(v) of the Theorem is proved.

\noindent$\bullet$ Let $u$ be a unit of order $11$. By (\ref{E:1})
and Proposition \ref{P:3} we get  $\nu_{11a}+\nu_{11b}=1$.
Applying Proposition \ref{P:1} to characters $\chi_{4}$ we obtain
\[
\begin{split}
\mu_{1}(u,\chi_{4},*) & = \textstyle \frac{1}{11} (6 \nu_{11a} -
5 \nu_{11b} + 16) \geq 0; \\ 
\mu_{2}(u,\chi_{4},*) & = \textstyle \frac{1}{11} (-5 \nu_{11a} +
6 \nu_{11b} + 16) \geq 0; \\ 
\mu_{1}(u,\chi_{4},3) & = \textstyle \frac{1}{11} (7 \nu_{11a} -
4 \nu_{11b} + 15) \geq 0; \\ 
\mu_{2}(u,\chi_{4},3) & = \textstyle \frac{1}{11} (-4 \nu_{11a} +
7 \nu_{11b} + 15) \geq 0, \\ 
\end{split}
\]
which admits only the four integer solutions listed in part (vi) of
the Theorem, such that  $\mu_{1}(u,\chi_{4},*)$,
$\mu_{2}(u,\chi_{4},*)$ and $\mu_{2}(u,\chi_{4},3)$ are non-negative
integers.

\noindent$\bullet$ Let $u$ be a unit of order $15$. By (\ref{E:1})
and Proposition \ref{P:3} we have that
$$
\nu_{3a}+\nu_{3b}+\nu_{5a}=1.
$$
Since $|u^5|=3$, for any character $\chi$ of $G$ we need to
consider five cases, defined by part (iii) of the Theorem. Put
\begin{equation}\label{E:3}
(\alpha,\beta,\gamma,\delta,\kappa) = \begin{cases}
(45,45,19,16,8), &\quad \text{if }\quad  \chi(u^5)=\chi(3a);\\
(51,42,13,22,11),&\quad \text{if }\quad  \chi(u^5)=\chi(3b); \\
(39,48,25,10,5),&\quad \text{if }\quad  \chi(u^5)=2\chi(3a)-\chi(3b);\\
(33,51,31,4,2),&\quad \text{if }\quad  \chi(u^5)=3\chi(3a)-2\chi(3b);\\
(57,39,7,28,14),&\quad \text{if }\quad  \chi(u^5)=-\chi(3a)+2\chi(3b).\\
\end{cases}
\end{equation}
By  (\ref{E:2}) we obtain the system of inequalities
$$
\mu_{0}(u,\chi_{6},*)= \textstyle \frac{1}{15}
(24\nu_{3b}+\alpha) \geq 0; \quad
\mu_{5}(u,\chi_{6},*)  = \textstyle\frac{1}{15}
(-12 \nu_{3b} + \beta) \geq 0.
$$
It follows  that  the integral solution is $\nu_{3b}\in\{-1,4\}$\;
if $(\alpha,\beta)=(39,48)$,  and
\begin{equation}\label{E:4}
\nu_{3b}=\begin{cases}
0, \quad &\text{if}\quad (\alpha,\beta)\in \{(45,45),(57,39)\};\\
1, \quad &\text{if}\quad (\alpha,\beta)=(51,42);\\
3,\quad &\text{if}\quad (\alpha,\beta)=(33,51).\\
\end{cases}
\end{equation}
\newline
Again by (2) we obtain the system of inequalities
\[
\begin{split}
\mu_{0}(u,\chi_{2},*) & = \textstyle \frac{1}{15} (16
\nu_{3a} - 8 \nu_{3b} +8 \nu_{5a} + \gamma) \geq 0; \\ 
\mu_{0}(u,\chi_{4},*) & = \textstyle \frac{1}{15} (-16 \nu_{3a} +
8 \nu_{3b} + 8 \nu_{5a} + \delta) \geq 0; \\ 
\mu_1(u,\chi_2,*) & = \textstyle \frac{1}{15} (2\nu_{3a} - \nu_{3b} + \nu_{5a} + \kappa)\geq 0,\\
\end{split}
\]
which has no integral solutions such that $(\alpha,\beta,\gamma,\delta,\kappa)$ from (\ref{E:3}),
$\nu_{3b}$ from (\ref{E:4}) and $\mu_{i}(u,\chi_{2},*)$ are non-negative integers.

$\bullet$ Let $u$ be a unit of  order $22$. By (\ref{E:1}) and
Proposition \ref{P:3} we have
$$
\nu_{2a}+\nu_{2b}+\nu_{11a}+\nu_{11b}=1.
$$
Since $u^{11}$ has order $2$ and $u^{2}$ has order $11$,  by (ii)
and (v) of the Theorem we have six and four different partial
augmentations for these orders, respectively. Thus we need to consider altogether $6\cdot
4=24$ cases. For any character $\chi$ of $G$ put
$$
\mathfrak{R}_\chi=\{\; \chi(11a),\quad  \chi(11b),\quad  2\chi(11a)-\chi(11b),\quad  -\chi(11a)+2\chi(11b)\; \};
$$
\[
(\alpha,\beta)=\begin{cases}
(10,12)&\quad \text{if }\;  \chi(u^2)\in \mathfrak{R}_\chi,\quad  \chi(u^{11})=\chi(2a);\\
(14,8)&\quad \text{if }\;  \chi(u^2)\in \mathfrak{R}_\chi,\quad  \chi(u^{11})=\chi(2b);\\
(22,0)&\quad \text{if }\;  \chi(u^2)\in \mathfrak{R}_\chi,\quad \chi(u^{11})=-2\chi(2a)+3\chi(2b);\\
(6,16)&\quad \text{if }\;  \chi(u^2)\in \mathfrak{R}_\chi,\quad \chi(u^{11})=2\chi(2a)-\chi(2b);\\
(2,20)&\quad \text{if }\;  \chi(u^2)\in \mathfrak{R}_\chi,\quad \chi(u^{11})=3\chi(2a)-2\chi(2b);\\
(18,4)&\quad \text{if }\;  \chi(u^2)\in \mathfrak{R}_\chi,\quad \chi(u^{11})=-\chi(2a)+2\chi(2b).\\
\end{cases}
\]
By (\ref{E:2}) we obtain that
\[
\begin{split}
\mu_0(u,\chi_2,*)&=\textstyle
\frac{1}{22}(-10\cdot (\nu_{2a}-3\nu_{2b})+\alpha)\geq 0; \\
\mu_{11}(u,\chi_2,*)&=\textstyle
\frac{1}{22}(10\cdot(\nu_{2a}-3\nu_{2b})+\beta)\geq 0.\\
\end{split}
\]
It is easy to check, that  if $(\alpha,\beta)\not\in\{(10,12),
(22,0)\}$, then  this system of inequalities has no integral
solutions such that $\mu_{i}(u,\chi_{2},*)$  are non-negative
integers.

If $(\alpha,\beta)=(10,12)$, then from  the last system of inequalities
we obtain that  $\nu_{2a}=3\nu_{2b}+1$, and put $(\gamma,\delta)=(50,60)$.
If $(\alpha,\beta)=(2,20)$, then from  the last system of inequalities
we get $\nu_{2a}=3\nu_{2b}-2$, and put $(\gamma,\delta)=(26,84)$.

Again, using (\ref{E:2})  in both cases of values of $(\alpha,\beta)$ we have
\[
\begin{split}
\mu_0(u,\chi_8,*) &= \textstyle \frac{1}{22} ( -50 \nu_{2a} + 70
\nu_{2b}+ \gamma )\geq  0;\\
\mu_{11}(u,\chi_8,*) &= \textstyle \frac{1}{22} ( 50
\nu_{2a}-70\nu_{2b}+ \delta )\geq 0,\\
\end{split}
\]
which has no integer solution such that $(\gamma,\delta)\in
\{(50,60),\; (26,84)\}$ and $\mu_{i}(u,\chi_{8},*)$ are non-negative
integers.

$\bullet$ Let $u$ be a unit of order $33$. Obviously, \quad
$\nu_{3a}+\nu_{3b}+\nu_{11a}+\nu_{11b}=1$\quad  by (\ref{E:1}) and
Proposition \ref{P:3}. Since $u^{11}$ has order $3$ and $u^3$ has
order $11$,  by (iii) and (v) of the Theorem we have five and four
different partial augmentations, respectively. Thus we need to
consider $20$ cases, such that
\[
\begin{split}
\chi(u^3)\in \mathfrak{R}_\chi=\{\;  -\chi(11a)+ 2\chi(11b),\quad
&\chi(11b),\\
& 2\chi(11a)-\chi(11b),\quad   \chi(11a)\; \};\\
\chi(u^{11})\in \{\; -\chi(3a)+ 2\chi(3b),\quad
&\chi(3a),\quad \chi(3b),\\
&3\chi(3a)-2\chi(3b),\quad 2\chi(3a)-\chi(3b)\; \},
\end{split}
\]
where $\chi$ is a character of the group $G$. Put
\begin{equation}\label{E:5}
(\alpha,\beta)=\begin{cases}
(15,9), &\quad \text{if } \chi(u^3)\in \mathfrak{R}_\chi,\quad \chi(u^{11})=\chi(3a);\\
(9,12),&\quad \text{if } \chi(u^3)\in \mathfrak{R}_\chi,\quad \chi(u^{11})=\chi(3b); \\
(21,6),&\quad \text{if } \chi(u^3)\in \mathfrak{R}_\chi,\quad \chi(u^{11})=2\chi(3a)-\chi(3b);\\
(27,3),&\quad \text{if } \chi(u^3)\in \mathfrak{R}_\chi,\quad \chi(u^{11})=3\chi(3a)-2\chi(3b);\\
(3,15),&\quad \text{if } \chi(u^3)\in \mathfrak{R}_\chi,\quad \chi(u^{11})=-\chi(3a)+2\chi(3b).\\
\end{cases}
\end{equation}
By (\ref{E:2})  we obtain  the system of inequalities
\[
\begin{split}
\mu_{0}(u,\chi_{2},*) & = \textstyle \frac{1}{33} (20\cdot ( 2\nu_{3a} -\nu_{3b}) +\alpha) \geq 0; \\ 
\mu_{11}(u,\chi_{2},*) & = \textstyle \frac{1}{33} (-10\cdot (2\nu_{3a}-\nu_{3b}) +\beta) \geq 0. \\ 
\end{split}
\]
It is easy to check that this system of inequalities has no
integral solutions such that $(\alpha,\beta)$ from (\ref{E:5}) and
$\mu_{i}(u,\chi_{2},*)$ are non-negative integers.

$\bullet$ Let $u$ be a unit of  order $55$. By (\ref{E:1}) and
Proposition \ref{P:3} we have that
$$
\nu_{5a}+\nu_{11a}+\nu_{11b}=1.
$$
Since $u^5$ has order $11$ and by (v) of the Theorem we have four
different partial augmentations, we need to consider the following
four cases:
\[
\begin{matrix}
\chi(u^5)&= \chi(11a);\qquad \chi(u^5)&=& -\chi(11a)+ 2\chi(11b);\\
\chi(u^5)&= \chi(11b);\qquad \chi(u^5)&=& \;2\chi(11a)-\chi(11b),
\end{matrix}
\]
where $\chi$ is a character of $G$. Again by (\ref{E:2})  we get  in all of
these four cases
\[
\begin{split}
\mu_{0}(u,\chi_{2},*)  & = \textstyle \frac{1}{55} (40 \nu_{5a} +15) \geq 0; \\
\mu_{1}(u,\chi_{2},*)  &= \textstyle \frac{1}{55} ( \nu_{5a}+10) \geq 0; \\ 
\mu_{11}(u,\chi_{2},*) & = \textstyle \frac{1}{55} (-10 \nu_{5a}+10) \geq 0. \\ 
\end{split}
\]
Clearly, this system of inequalities has no integral solutions such
that $\mu_{0}(u,\chi_{2},*)$, $\mu_{1}(u,\chi_{2},*)$ and
$\mu_{11}(u,\chi_{2},*)$ are non-negative integers. The proof is
done.

\qed

\subsection*{Acknowledgment}
The authors are grateful to the referee for his useful comments.

\bibliographystyle{plain}
\bibliography{Bovdi_Konovalov_Siciliano_M12}

\end{document}